\documentclass[10pt]{amsart}

\usepackage{latexsym,amssymb,amsmath,graphics}
\usepackage[dvips]{graphicx}

\def\ff{{\mathcal F}}

\def\ffi{\varphi}

\def\dst{\displaystyle}

\def\R{{\mathbb{R}}}

\def\Z{{\mathbb{Z}}}

\newcommand{\norm}[1]{{\left\|{#1}\right\|}}

\newcommand{\abs}[1]{{\left|{#1}\right|}}

\newenvironment{notation}[1][]{\vskip1pt\noindent\rm\textit{Notation}\,:\ }{\rm\vskip1pt}
\newenvironment{remark}[1][]{\vskip1pt\noindent\rm\textit{Remark #1}\,:\ }{\rm\vskip1pt}
\newenvironment{definition}[1][]{\vskip3pt\noindent\sl\textbf{Definition.}\ }{\rm\vskip3pt}

\newtheorem{lemma}{Lemma}[section]

\newtheorem{theorem}[lemma]{Theorem}
\newtheorem{corollary}[lemma]{Corollary}

\begin{document}

\title{Zero-free regions of radar ambiguity functions and moments}
\author{Philippe Jaming}
\address{Universit\'e d'Orl\'eans\\
Facult\'e des Sciences\\ 
MAPMO - F\'ed\'eration Denis Poisson\\ BP 6759\\ F 45067 Orl\'eans Cedex 2\\
France}
\email{Philippe.Jaming@univ-orleans.fr}

\begin{abstract}
In this article, we give an estimate of the zero-free region around the origin of the ambiguity function
of a signal $u$ in terms of the moments of $u$. This is done by proving an uncertainty relation between the first zero of the Fourier
transform of a non-negative function and the moments of the function. As a corollary, we also give an estimate of how much
a function needs to be translated to obtaina function that is orthogonal to the original function.
\end{abstract}




\thanks{Research partially financed by: {\it European Commission}
Harmonic Analysis and Related Problems 2002-2006 IHP Network
(Contract Number: HPRN-CT-2001-00273 - HARP)
}
\maketitle

\section{Introduction}\ \\

Woodward's time-frequency correlation function or radar ambiguity function \cite{ieee.Wo,ieee.Wi,ieee.AT},
as defined by
$$
A(u)(x,y)=\int_{-\infty}^{+\infty}u\left(t+\frac{x}{2}\right)\overline{u\left(t-\frac{x}{2}\right)}
e^{-2i\pi yt}\mbox{d}t
$$
plays a central role in evaluating the ability of a transmitted radar waveform $\mbox{Re}\,\bigl(u(t)e^{i\omega_0 t}\bigr)$
to distinguish targets that are separated by range delay $x$ and Doppler frequency $y$.
Ideally, one would like $A(u)$ to be a Dirac mass at $(0,0)$, but this desideratum is not achievable because 
of the ``\emph{ambiguity uncertainty principle}'', that is the constraint
$$
\iint_{\R^2}|A(u)(x,y)|^2\mbox{d}x\,\mbox{d}y=A(u)(0,0)^2=\left(\int_\R|u(t)|^2\mbox{d}t\right)^2.
$$
As $A(u)$ is continuous when $u$ is a signal of finite energy, it follows that $A(u)$ can not vanish in a neighborhood
of $(0,0)$. 
Since ideal behavior is not achievable, it becomes important to determine how closely to the ideal situation
one can come. 

A major attempt in that direction is due to Price and Hofstetter \cite{ieee.PH} who considered
the quantity
$$
V(E)=\iint_E|A(u)(x,y)|^2\mbox{d}x\,\mbox{d}y
$$
where $E$ is a measurable subset of $\R^2$. It should however be mentioned that, for their results
to be significant, one has to go outside the class
of signals of finite energy as $V(E)$ is supposed to have a limit $V_0$ when $E$ shrinks to $\{(0,0)\}$.
Indeed, when $u$ has finite energy, from the continuity of $A(u)$, one gets that 
$$
V(E)\leq \mbox{area}(E)\sup_{(x,y)\in E}|A(u)(x,y)|^2\to \mbox{area}(\{(0,0)\})|A(u)(0,0)|^2=0.
$$
Nevertheless, it is possible to define $A(u)$ when $u$ is a Schwartz distribution and the assumption
in \cite[Section II]{ieee.PH} is that $A(u)\in L^2_{\mathrm{loc}}(\R^2\setminus\{(0,0)\})$.

In this paper, we restrict our attention to signals $u$ of finite energy and we try to determine more precisely
the neighborhood of $(0,0)$ on which $A(u)$ does not vanish. To do so, we
prove a new form of uncertainty principle,
showing that there is an exclusion relation between the function $u$ having moments and its ambiguity function
being $0$ near $(0,0)$. This result is inspired by a recent result of Luo and Zhang
concerning the Fourier transform. They prove that if a real non-negative valued function is supported in $[0,+\infty)$
then there is an uncertainty principle relating its moments and the first zero its Fourier transform.

It turns out that the ambiguity function, when restricted in a given direction, is always the Fourier
transform of a non-negative function, but the support condition of Luo and Zhang is not valid.
We thus start by removing that condition in their uncertainty principle, which can be done at little expenses
apart from some numerical constants. This then allows to obtain a zero-free region for the ambiguity function
$A(u)$ when $u$ and its Fourier transforms have $L^2$-moments.
This region turns out to be a square when one considers dispersions.

The article is divided in two sections. The first one is devoted to the extension of Luo and Zhang's
uncertainty principle. In the second section, we recall how the fractional Fourier transform allows to see
the restriction of the ambiguity function in a given direction as a Fourier transform of a non-negative 
function. We then apply the results of the first section to obtain zero-free regions of the ambiguity function.

\section{Zero-free regions of the Fourier transform}
\begin{notation} 
For $1\leq p<\infty$ we define $L^p(\R)$ as the space of measurable functions such that
$$
\norm{f}_p^p:=\int_\R|f(t)|^p\,\mbox{d}t<+\infty.
$$
For $u\in L^1(\R)\cap L^2(\R)$ we define the Fourier transform as
$$
\ff u(\xi)=\widehat{u}(\xi)=\int_{\R}u(t)e^{-2i\pi\xi t}\,\mbox{d}t,\quad\xi\in\R
$$
and then extend it to $L^2(\R)$ in the usual way.
\end{notation}

\begin{theorem}
\label{ieee:th1}\ \\
For every $q>0$ there exists $\kappa_q>0$ such that,
if $u\in L^1(\R)$ is a non-negative function such that $u\not=0$, then
$$
\inf\{\xi>0\,:\ \widehat{u}(\pm\xi)=0\}^q\inf_{t_0\in\R}\bigl\||t-t_0|^qu(t)\bigr\|_1\geq\kappa_q\norm{u}_1.
$$
\end{theorem}

\begin{remark}\ \\
--- One may take for the constant $\kappa=\frac{1}{c(2\pi)^q}$ where $c$ is the smallest constant
in Equation \eqref{ieee:eq3} below.\\
--- A similar result has been proved in \cite{ieee.LZ} but with the extra assumption
that $u$ be supported in $[0,+\infty)$. The constant $\kappa$ is then explicitely
known and better than the one above.
\end{remark}

\begin{proof}[Proof of Theorem \ref{ieee:th1}]
The proof is very similar to that of Luo and Zhang. The idea is that
near $0$, $e^{ix}\sim\cos x\sim 1-x^2/2$ and that for $x$ big enough
$1-x^2/2\ll\cos x$. Indeed $\mbox{Re}\,e^{ix}=\cos x\geq 1-x^2/2$ for all $x$.
Thus $\dst\mbox{Re}\,\widehat{u}\geq\int_\R\bigl(1-(2\pi t\xi)^2/2\bigr)u(t)\,\mbox{d}t$.
It follows that if $\widehat{u}(\xi)=0$, then $\mbox{Re}\,\widehat{u}=0$ thus
$$
\xi^2\norm{t^2u(t)}_1\geq\frac{1}{2\pi^2}\norm{u(t)}_1.
$$
The constant $\frac{1}{2\pi^2}$ can be improved when $u$ is supported in $[0,+\infty)$
using a refined version of the inequality $\cos x\geq 1-x^2/2$
({\it see} \cite[Proposition 1.2]{ieee.LZ})
and can further be extended by replacing $\norm{t^2u}$ by other moments $\norm{t^qu}$.
Our substitute to \cite[Proposition 1.2]{ieee.LZ} is the following~:

\medskip

\noindent{\bf Fact.} {\sl for every $q>0$, there exists $a,c$ such that, for all $x\in\R$,}
\begin{equation}
\label{ieee:eq3}
a\cos x\geq 1-c|x|^q.
\end{equation}

\medskip

It then follows that, if $\xi$ is such that $\widehat{u}(\xi)=0$, then
\begin{eqnarray*}
0&=&a\mbox{Re}\,\widehat{u}(\xi)
=\int u(t)a\cos 2\pi t\xi\,\mbox{d}t
\geq\int u(t)\Bigl(1-c(2\pi|\xi|)^q|t|^q\Bigr)\,\mbox{d}t\\
&=&\norm{u}_1-c(2\pi|\xi|)^q\|\,|t|^qu\,\|_1.
\end{eqnarray*}

Now, as $\widehat{u}$ is continuous, we may take $\xi=\pm \tau$, and get that
$$
|\tau|^q\,\bigl\|\,|t|^qu\,\bigr\|_1\geq\frac{1}{c(2\pi)^q}\bigl\|u\bigr\|_1
$$
and applying this to the translate $u_{t_0}(t)=u(t-t_0)$ we get the desired result.
\end{proof}

\begin{proof}[Proof of \eqref{ieee:eq3}]
The fact is trivial as long as one does not look for best constants. Indeed, one may take $a=2$ and $c$ such that $1-c\abs{\frac{\pi}{3}}^q=-2$
{\it i.e.} $c=3\left(\frac{3}{\pi}\right)^q$. In this case, for $\dst 0\leq|x|\leq\frac{\pi}{3}$,
$1-c\abs{x}^q\leq 1\leq 2\cos x$ and for $|x|>\frac{\pi}{3}$, $1-c\abs{x}^q\leq -2\leq 2\cos x$.
A slightly more refined argument, taking $a=1+\eta$, $x_0=\arccos(1+\eta)^{-1}$ gives
$c=\dst\frac{2+\eta}{\bigl(\arccos(1+\eta)^{-1}\bigr)^q}$. One may then minimize over $\eta>0$
and some values are given in the following table~:
\begin{center}
\begin{tabular}{c|cccc}
q&3&4&5&6\\
\hline
a&3.26&3.94&4.61&5.27\\
\hline
c&2.134&1.656&1.241&0.908\\
\end{tabular}
\end{center}
The constant $c$ may be slightly improved
by a more refined argument. More precisely, at the point $x_0$ above, $1-c|x|^q$ and $a\cos x$
are still far apart. Optimizing between all parameters is however difficult.
We would like to mention that, even for $q=2$, the estimate $\cos x\geq 1-x^2/2$ is not best possible
for our needs. Indeed, choosing $a=1$, which is the smallest possible value for $a$ when $q\leq 2$, 
may not lead to the best choice of $c$.
For instance, a computer plot will convince the reader that
$$
1.02\cos x\geq 1-0.52x^{3/2}\quad,\quad
1.1\cos x\geq 1-0.42x^2.
$$
However, for $q\leq 1$, $a=1$ allows for the best constant $c$ since $1-cx^q$ is then concave on $]0,+\infty)$.
The best constant $c$ can then
be computed as follows~: the equation $\cos x=1-cx^q$ has to have a solution $x_0$
in $[\pi/2,\pi]$ for which $\sin x_0=cq|x_0|^{q-1}$. It then follows that $x_0$ is the unique solution
in $[\pi/2,\pi]$ of
$$
\cos x+\frac{1}{q}x\sin x=1\quad\mbox{and then }c=\frac{\sin x_0}{q x_0^{q-1}}.
$$
This equation has a solution since for $x=\pi/2$, $\cos x+\frac{1}{q}x\sin x=\frac{\pi}{2q}>1$
and for $x=\pi$, $\cos x+\frac{1}{q}x\sin x=-1<1$ and $\ffi(x)=\cos x+\frac{1}{q}x\sin x$ is continuous.
The solution is unique since $\ffi'(x)=\left(\frac{1}{q}-1\right)\sin x+\frac{1}{q}x\cos x$
is made of 2 pieces $\ffi_1(x)=\left(\frac{1}{q}-1\right)\sin x$ and $\ffi_2(x)=\frac{1}{q}x\cos x$.
The first one, $\ffi_1$ is non-negative and decreasing on $[\pi/2,\pi]$ while the second one
is negative decreasing on $[\pi/2,\pi]$. As $\ffi'(\pi/2)\geq0$ and $\ffi'(\pi)\leq0$, there exists $x_1$ such that
$\ffi'(x_1)=0$ and $\ffi'(x)>0$ for $\pi/2<x<x_1$ while
$\ffi'(x)<0$ for $x_1<x<\pi$. It follows that the solution of $\ffi(x)=0$ in $[x_1,\pi]$ is unique.

Finally note also that for $q\sim 0$, $c\sim 2$ and for $q\sim 1$, $c\sim 0.73$.
\end{proof}

To conclude this section, let us give a first application of this result. We ask whether a translate $f_a(t)=f(t-a)$  of $f$
can be orthogonal to $f$. But
$$
0=\int_\R f(t)\overline{f(t-a)}\,\mbox{d}t=\int_\R |\widehat{f}(\xi)|^2e^{2ia\xi}\,\mbox{d}\xi=\ff[|\widehat{f}|^2](-a).
$$
Similarily, if the modulation $f^{(\omega)}(t)=e^{2i\pi\omega t}f(t)$ of $f$ is orthogonal to $f$ then $\ff[|f|^2](\omega)=0$
From Theorem \ref{ieee:th1}, we get the following:

\begin{corollary}\ \\
Let $q>0$ and $\kappa_q$ be the constant of Theorem \ref{ieee:th1}. Let $f\in L^2(\R)$.

\noindent--- Assume that $(1+|\xi|)^{q/2}\widehat{f}\in L^2$. Then for $f$ and its translate $f_a$ to be orthogonal, it is necessary that
$$
|a|^q \inf_{t_0\in\R}\|\,|t-t_0|^{q/2}\widehat{f}(t)\|_2^2\geq \kappa_q\norm{f}_2^2.
$$
--- Assume that $(1+|t|)^{q/2}f\in L^2$. Then for $f$ and its modulation $f^{(\omega)}$ to be orthogonal, it is necessary that
$$
|\omega|^q \inf_{t_0\in\R}\|\,|t-t_0|^{q/2}f(t)\|_2^2\geq \kappa_q\norm{f}_2^2.
$$
\end{corollary}

\section{Zero-free regions of the ambiguity function}

\subsection{Fractional Fourier transforms}\ \\[3pt]
For $\alpha\in\R\setminus\pi\Z$, let $c_\alpha=\dst\frac{\exp\frac{i}{2}\left(\alpha-\frac{\pi}{4}\right)}{\sqrt{|\sin\alpha|}}$
be a square root of $1-i\cot\alpha$. For $f\in L^1(\R)$ and $\alpha\notin\pi\Z$, define
$$
\ff_\alpha f(\xi)=c_\alpha e^{-i\pi\xi^2\cot\alpha}\int_\R f(t)e^{-i\pi t^2\cot\alpha}e^{-2i\pi t\xi/\sin\alpha}\mbox{d}t
=c_\alpha e^{-i\pi\xi^2\cot\alpha}\ff[f(t)e^{-i\pi t^2\cot\alpha}](\xi/\sin\alpha)
$$
while for $k\in\Z$, $\ff_{2k\pi} f=f$ and $\ff_{(2k+1)\pi}f(\xi)=f(-\xi)$.
This transformation has the following properties~:
\begin{enumerate}
\item\label{ieee:prop:fa1} $\dst\int_\R\ff_\alpha f(\xi)\overline{\ff_\alpha g(\xi)}\mbox{d}\xi=\int_\R f(t)\overline{g(t)}\mbox{d}t$
which allows to extend $\ff_\alpha$ from $L^1(\R)\cap L^2(\R)$ to
$L^2(\R)$ as a unitary operator on $L^2(\R)$;

\item\label{ieee:prop:fa2} $\ff_\alpha\ff_\beta=\ff_{\alpha+\beta}$;

\item\label{ieee:prop:fa3} if $f_a(t)=f(t-a)$ then
$$
\ff_\alpha f_a(\xi)=\ff_\alpha f(\xi+a\cos\alpha)e^{-i\pi a^2\cos\alpha\sin\alpha-2i\pi a\xi\sin\alpha};
$$
\item\label{ieee:prop:fa4} if $f_\omega(t)=e^{-2i\pi\omega t}f(t)$ then
$$
\ff_\alpha f_\omega(\xi)=\ff_\alpha f(\xi+\omega\sin\alpha)e^{i\pi\omega^2\cos\alpha\sin\alpha+2i\pi\omega\xi\sin\alpha};
$$
\item\label{ieee:prop:fa5} if $f\in L^2$ is such that $tf\in L^2$ then 
$$
\ff_\alpha[tf](\xi)=\xi\ff_\alpha f(\xi)\cos\alpha+i[\ff_\alpha f]'(\xi)\sin\alpha.
$$
\end{enumerate}

Let us recall that the ambiguity function of $u\in L^2(\R)$ is defined by
$$
A(u)(x,y)=\int_\R u\left(t+\frac{x}{2}\right)\overline{u\left(t-\frac{x}{2}\right)}
e^{-2i\pi ty}\mbox{d}t.
$$
The following properties are well known \cite{ieee.Al,ieee.AT,ieee.Wi}~:
\begin{enumerate}
\item\label{ieee:prop:amb1} $A(u)\in L^2(\R^2)$ with $\norm{A(u)}_{L^2(\R^2)}=\norm{u}_2^2$ and is continuous;

\item\label{ieee:prop:amb2} $A(u)(0,0)=\norm{u}_2^2$ where it is maximal;

\item\label{ieee:prop:amb3} $A(u)(-x,-y)=\overline{A(u)(x,y)}$;

\item\label{ieee:prop:amb4} $A(\ff_\alpha u)(x,y)=A(u)(x\cos\alpha-y\sin\alpha,x\sin\alpha+y\cos\alpha)$.
\end{enumerate}

The second property was proved in \cite{ieee.Wi} when the fractional Fourier transform is defined in terms of Hermite polynomials
and in \cite{ieee.Al} with the above definition of the fractional Fourier transform.

\subsection{Zero free regions}\ \\
Noticing that $A(u)$ is a Fourier transform and in particular that $A(u)(0,y)=\ff[|u|^2](y)$,
we get from Property \ref{ieee:prop:amb4} of the ambiguity function that
$$
A(u)(-y\sin\alpha,y\cos\alpha)=\ff[|\ff_\alpha u|^2](y).
$$
\begin{definition}\ \\
Let us define, for $\theta\in]0,\pi[$
$$
\tau_\theta=\inf\{t>0~:A(u)(t\cos\theta,t\sin\theta)=0\mbox{ or }A(u)(-t\cos\theta,-t\sin\theta)=0\}.
$$
\end{definition}\ \\
Let $q>0$ and $\kappa_q$ be given by Theorem \ref{ieee:th1}.
Applying this theorem, we obtain
\begin{equation}
\label{ieee.eq.tau1}
\tau_\theta^q\inf\limits_{t_0\in\R}\norm{|t-t_0|^q|\ff_{\theta-\pi/2} u|^2}_1\geq\kappa_q
\norm{|\ff_{\theta-\pi/2} u|^2}_1=\kappa_q\norm{u}_2^2.
\end{equation}

Let us now show that in the specific case, a more precise result can be obtained:
if $u\in L^2$ is such that $tu\in L^2$ and $t\widehat{u}\in L^2$
then from Property \ref{ieee:prop:fa5} of the fractional Fourier transform
\begin{eqnarray*}
\norm{t\ff_\alpha u}_2&=&\norm{tu(t)\cos(-\alpha)-iu'\sin\alpha}_2
\leq\norm{tu(t)}_2|\cos\alpha|+\norm{u'}_2|\sin\alpha|\\
&=&\norm{tu(t)}_2|\cos\alpha|+\norm{\xi\widehat{u}(\xi)}_2|\sin\alpha|.
\end{eqnarray*}
In particular, the ambiguity function $A(u)$ of $u$ has no zero in the region
$$
\left\{(t\sin\alpha,-t\cos\alpha)\,: 0<\alpha<\pi,\ 
|t|\leq\frac{\sqrt{2}\norm{u}_2}{2\pi\bigl(\norm{tu(t)}_2|\cos\alpha|+\norm{\xi\widehat{u}(\xi)}_2|\sin\alpha|\bigr)}
\right\}.
$$
This region is a rombus with endpoints
$$
\left(\pm\frac{\sqrt{2}\norm{u}_2}{2\pi\bigl\|tu(t)\bigr\|_2},0\right)\quad\mbox{and}\quad
\left(0,\pm\frac{\sqrt{2}\norm{u}_2}{2\pi\bigl\|\xi\widehat{u}(t)\bigr\|_2}\right).
$$
Further, changing $u(t)$ into $u(t-a)e^{i\omega t}$ leaves the modulus of $A(u)$ unchanged,
and so are the zero-free regions of $A(u)$. We have thus proved

\begin{theorem}\ \\
Let $u\in L^2(\R)$ be such that $tu(t)\in L^2(\R)$ and $\xi\widehat{u}(\xi)\in L^2(\R)$.
Then the ambiguity function $A(u)$ of $u$ has no zero in the convex hull of the four points
$$
\left(\pm\frac{\sqrt{2}\norm{u}_2}{2\pi\inf_{a\in\R}\bigl\||t-a|u(t)\bigr\|_2},0\right)\quad\mbox{and}\quad
\left(0,\pm\frac{\sqrt{2}\norm{u}_2}{2\pi\inf_{\omega\in\R}\bigl\||\xi-\omega|\widehat{u}(t)\bigr\|_2}\right).
$$
\end{theorem}

The area of that rombus is
$$
\frac{\frac{1}{\pi}\norm{u}_2^2}{\inf_{a\in\R}\bigl\||t-a|u(t)\bigr\|_2\inf_{\omega\in\R}\bigl\||\xi-\omega|\widehat{u}(t)\bigr\|_2}
\leq 4
$$
according to Heisenberg's uncertainty principle.

Note also that the numerical constant $\frac{\sqrt{2}}{2\pi}$ can be improved to $0.248$ using
the inequality $1.1\cos x\geq 1-0.41x^2$ instead of $\cos x\geq 1-x^2/2$.

\end{document}